\documentclass[12pt]{article}

\usepackage{amsmath,amssymb,amsthm,amscd,a4wide,upref}

\usepackage[enableskew]{youngtab}

\usepackage{hyperref}
\hypersetup{colorlinks,citecolor=blue,filecolor=black,linkcolor=blue,urlcolor=blue}
\usepackage{enumerate}
\usepackage{mathtools}
\usepackage{enumitem}

\usepackage{lmodern}     % set math font to Latin modern math
\usepackage[T1]{fontenc}
\begin{document}

% % % % % % % % % % % % % %
%\newcommand{\tma}{\textcolor{magenta}}
%\newcommand{\tre}{\textcolor{red}}
% % % % % % % % % % % % % %

\newcommand{\ad}{{\rm ad}}
\newcommand{\cri}{{\rm cri}}
\newcommand{\row}{{\rm row}}
\newcommand{\col}{{\rm col}}
\newcommand{\End}{{\rm{End}\ts}}
\newcommand{\Rep}{{\rm{Rep}\ts}}
\newcommand{\Hom}{{\rm{Hom}}}
\newcommand{\Mat}{{\rm{Mat}}}
\newcommand{\ch}{{\rm{ch}\ts}}
\newcommand{\chara}{{\rm{char}\ts}}
\newcommand{\diag}{{\rm diag}}
\newcommand{\st}{{\rm st}}
\newcommand{\non}{\nonumber}
\newcommand{\wt}{\widetilde}
\newcommand{\wh}{\widehat}
\newcommand{\ol}{\overline}
\newcommand{\ot}{\otimes}
\newcommand{\la}{\lambda}
\newcommand{\La}{\Lambda}
\newcommand{\De}{\Delta}
\newcommand{\al}{\alpha}
\newcommand{\be}{\beta}
\newcommand{\ga}{\gamma}
\newcommand{\Ga}{\Gamma}
\newcommand{\ep}{\epsilon}
\newcommand{\ka}{\kappa}
\newcommand{\vk}{\varkappa}
\newcommand{\si}{\sigma}
\newcommand{\vs}{\varsigma}
\newcommand{\vp}{\varphi}
\newcommand{\de}{\delta}
\newcommand{\ze}{\zeta}
\newcommand{\om}{\omega}
\newcommand{\Om}{\Omega}
\newcommand{\ee}{\epsilon^{}}
\newcommand{\su}{s^{}}
\newcommand{\hra}{\hookrightarrow}
\newcommand{\ve}{\varepsilon}
\newcommand{\ts}{\,}
\newcommand{\vac}{\mathbf{1}}
\newcommand{\vacu}{|0\rangle}
\newcommand{\di}{\partial}
\newcommand{\qin}{q^{-1}}
\newcommand{\tss}{\hspace{1pt}}
\newcommand{\Sr}{ {\rm S}}
\newcommand{\U}{ {\rm U}}
\newcommand{\BL}{ {\overline L}}
\newcommand{\BE}{ {\overline E}}
\newcommand{\BP}{ {\overline P}}
\newcommand{\AAb}{\mathbb{A}\tss}
\newcommand{\CC}{\mathbb{C}\tss}
\newcommand{\KK}{\mathbb{K}\tss}
\newcommand{\QQ}{\mathbb{Q}\tss}
\newcommand{\SSb}{\mathbb{S}\tss}
\newcommand{\TT}{\mathbb{T}\tss}
\newcommand{\ZZ}{\mathbb{Z}\tss}
\newcommand{\DY}{ {\rm DY}}
\newcommand{\X}{ {\rm X}}
\newcommand{\Y}{ {\rm Y}}
\newcommand{\Z}{{\rm Z}}
\newcommand{\Ac}{\mathcal{A}}
\newcommand{\Lc}{\mathcal{L}}
\newcommand{\Mc}{\mathcal{M}}
\newcommand{\Pc}{\mathcal{P}}
\newcommand{\Qc}{\mathcal{Q}}
\newcommand{\Rc}{\mathcal{R}}
\newcommand{\Sc}{\mathcal{S}}
\newcommand{\Tc}{\mathcal{T}}
\newcommand{\Bc}{\mathcal{B}}
\newcommand{\Ec}{\mathcal{E}}
\newcommand{\Fc}{\mathcal{F}}
\newcommand{\Gc}{\mathcal{G}}
\newcommand{\Hc}{\mathcal{H}}
\newcommand{\Uc}{\mathcal{U}}
\newcommand{\Vc}{\mathcal{V}}
\newcommand{\Wc}{\mathcal{W}}
\newcommand{\Yc}{\mathcal{Y}}
\newcommand{\Cl}{\mathcal{C}l}
\newcommand{\Ar}{{\rm A}}
\newcommand{\Br}{{\rm B}}
\newcommand{\Ir}{{\rm I}}
\newcommand{\Fr}{{\rm F}}
\newcommand{\Jr}{{\rm J}}
\newcommand{\Or}{{\rm O}}
\newcommand{\GL}{{\rm GL}}
\newcommand{\Spr}{{\rm Sp}}
\newcommand{\Rr}{{\rm R}}
\newcommand{\Zr}{{\rm Z}}
\newcommand{\gl}{\mathfrak{gl}}
\newcommand{\middd}{{\rm mid}}
\newcommand{\ev}{{\rm ev}}
\newcommand{\Pf}{{\rm Pf}}
\newcommand{\Norm}{{\rm Norm\tss}}
\newcommand{\oa}{\mathfrak{o}}
\newcommand{\spa}{\mathfrak{sp}}
\newcommand{\osp}{\mathfrak{osp}}
\newcommand{\f}{\mathfrak{f}}
\newcommand{\g}{\mathfrak{g}}
\newcommand{\h}{\mathfrak h}
\newcommand{\n}{\mathfrak n}
\newcommand{\m}{\mathfrak m}
\newcommand{\z}{\mathfrak{z}}
\newcommand{\Zgot}{\mathfrak{Z}}
\newcommand{\p}{\mathfrak{p}}
\newcommand{\sll}{\mathfrak{sl}}
\newcommand{\agot}{\mathfrak{a}}
\newcommand{\bgot}{\mathfrak{b}}
\newcommand{\qdet}{ {\rm qdet}\ts}
\newcommand{\Ber}{ {\rm Ber}\ts}
\newcommand{\HC}{ {\mathcal HC}}
\newcommand{\cdet}{{\rm cdet}}
\newcommand{\rdet}{{\rm rdet}}
\newcommand{\tr}{ {\rm tr}}
\newcommand{\gr}{ {\rm gr}\ts}
\newcommand{\str}{ {\rm str}}
\newcommand{\loc}{{\rm loc}}
\newcommand{\Gr}{{\rm G}}
\newcommand{\sgn}{ {\rm sgn}\ts}
\newcommand{\sign}{{\rm sgn}}
\newcommand{\ba}{\bar{a}}
\newcommand{\bb}{\bar{b}}
\newcommand{\bi}{\bar{\imath}}
\newcommand{\bj}{\bar{\jmath}}
\newcommand{\bk}{\bar{k}}
\newcommand{\bl}{\bar{l}}
\newcommand{\hb}{\mathbf{h}}
\newcommand{\Sym}{\mathfrak S}
\newcommand{\fand}{\quad\text{and}\quad}
\newcommand{\Fand}{\qquad\text{and}\qquad}
\newcommand{\For}{\qquad\text{or}\qquad}
\newcommand{\for}{\quad\text{or}\quad}
\newcommand{\grpr}{{\rm gr}^{\tss\prime}\ts}
\newcommand{\degpr}{{\rm deg}^{\tss\prime}\tss}
\newcommand{\bideg}{{\rm bideg}\ts}

\renewcommand{\theequation}{\arabic{section}.\arabic{equation}}

\numberwithin{equation}{section}

\newtheorem{thm}{Theorem}[section]
\newtheorem{lem}[thm]{Lemma}
\newtheorem{prop}[thm]{Proposition}
\newtheorem{cor}[thm]{Corollary}
\newtheorem{conj}[thm]{Conjecture}
\newtheorem*{mthm}{Main Theorem}
\newtheorem*{mthma}{Theorem A}
\newtheorem*{mthmb}{Theorem B}
\newtheorem*{mthmc}{Theorem C}
\newtheorem*{mthmd}{Theorem D}

\theoremstyle{definition}
\newtheorem{defin}[thm]{Definition}

\theoremstyle{remark}
\newtheorem{remark}[thm]{Remark}
\newtheorem{example}[thm]{Example}
\newtheorem{examples}[thm]{Examples}

\newcommand{\bth}{\begin{thm}}
\renewcommand{\eth}{\end{thm}}
\newcommand{\bpr}{\begin{prop}}
\newcommand{\epr}{\end{prop}}
\newcommand{\ble}{\begin{lem}}
\newcommand{\ele}{\end{lem}}
\newcommand{\bco}{\begin{cor}}
\newcommand{\eco}{\end{cor}}
\newcommand{\bde}{\begin{defin}}
\newcommand{\ede}{\end{defin}}
\newcommand{\bex}{\begin{example}}
\newcommand{\eex}{\end{example}}
\newcommand{\bes}{\begin{examples}}
\newcommand{\ees}{\end{examples}}
\newcommand{\bre}{\begin{remark}}
\newcommand{\ere}{\end{remark}}
\newcommand{\bcj}{\begin{conj}}
\newcommand{\ecj}{\end{conj}}

\newcommand{\bal}{\begin{aligned}}
\newcommand{\eal}{\end{aligned}}
\newcommand{\beq}{\begin{equation}}
\newcommand{\eeq}{\end{equation}}
\newcommand{\ben}{\begin{equation*}}
\newcommand{\een}{\end{equation*}}

\newcommand{\bpf}{\begin{proof}}
\newcommand{\epf}{\end{proof}}

\def\beql#1{\begin{equation}\label{#1}}

\newcommand{\Res}{\mathop{\mathrm{Res}}}

\title{\Large\bf $\Wc$-algebras associated with centralizers in type $A$}

\author{A. I. Molev}

\date{} % Start January 2020
\maketitle

%\vspace{4 mm}

\begin{abstract}
We introduce a new family of affine $\Wc$-algebras $\Wc^{\tss k}(\agot)$ associated with
the centralizers of arbitrary nilpotent elements in $\gl_N$.
We define them by using a version of the BRST complex of the quantum
Drinfeld--Sokolov reduction.
A family of free generators of $\Wc^{\tss k}(\agot)$ is produced in an explicit form.
We also give an analogue of the Fateev--Lukyanov realization
for the new $\Wc$-algebras by applying
a Miura-type map.
\end{abstract}

%\vspace{5 mm}
%%%
%%%{\it Key words:}
%%%

%\newpage

%\tableofcontents
%
%\newpage

\section{Introduction}
\label{sec:int}

The {\em affine
$\Wc$-algebra $\Wc^{\tss k}(\g)$ at the level $k\in\CC$}
associated with a simple Lie algebra $\g$
is a vertex algebra defined by
a quantum Drinfeld--Sokolov reduction
\cite{ff:qd}. These algebras originate in conformal field theory
and first appeared
in the work of Zamolodchikov~\cite{z:ie} and Fateev and Lukyanov~\cite{fl:mt}.
They were intensively studied both in mathematics and physics literature; see e.g.
\cite{a:rtw}, \cite{a:iw}, \cite{bs:ws}, \cite[Ch.~15]{fb:va} for detailed reviews.
More general $\Wc$-algebras $\Wc^{\tss k}(\g,f)$
were introduced in \cite{krw:qr}, which depend on a simple Lie (super)algebra $\g$ and
an (even) nilpotent element $f\in \g$ so that $\Wc^{\tss k}(\g)$ corresponds to
a principal nilpotent element $f$. Their counterparts for odd nilpotents
$f$ were studied in \cite{mr:qh} and \cite{mrs:sw} from the viewpoint of
quantum hamiltonian reduction.

In an important particular case, where the level $k$ takes the {\em critical value}, the
vertex algebra $\Wc^{\ts\cri}(\g)$ is commutative. It is isomorphic to
the center $\z(\wh\g)$ of the affine vertex algebra $V^{\cri}(\g)$ and is known
as the {\em Feigin--Frenkel center}
following the paper \cite{ff:ak}, where the structure of $\z(\wh\g)$ was described;
see also \cite{f:lc} for detailed arguments and \cite{m:so} for explicit constructions
of generators of $\z(\wh\g)$ for the Lie algebras $\g$ of classical types.
The arguments rely on a basic property of the simple Lie algebras stating that the
subalgebra of $\g$-invariants of the symmetric algebra $\Sr(\g)$ is a polynomial algebra.
It was shown in \cite{ppy:si}, that this property
is shared by a wide class of non-reductive Lie algebras $\agot$ which are centralizers
of nilpotent elements in $\g$. This served as a starting point for the work \cite{ap:qm},
where the Feigin--Frenkel theorem was extended to the affine vertex algebras
$V^{\cri}(\agot)$ associated with the centralizers.

These results motivate the question whether analogues of the $\Wc$-algebras
can be associated with the underlying Lie algebras $\agot$.
Our goal in this paper is to introduce and describe some basic properties
of $\Wc$-algebras $\Wc^{\tss k}(\agot)$, where $\agot$
is the centralizer of a nilpotent element $e$ in $\gl_N$. In the case $e=0$
the corresponding algebra coincides with the principal $\Wc$-algebra $\Wc^{\tss k}(\gl_N)$.

Since the Lie algebra $\agot$ is not semisimple,
the construction depends on a choice of the invariant symmetric
bilinear form on $\agot$. We will follow \cite{ap:qm}, where a natural form was used
to introduce the corresponding affine Kac--Moody algebra $\wh\agot$.
Its vacuum module $V^k(\agot)$ at the level $k$ is a vertex algebra.
The Lie algebra $\agot$ admits a triangular decomposition $\agot=\n_-\oplus\h\oplus \n_+$
which gives rise to a Clifford algebra
associated with $\n_+$ and we let $\Fc$ be its vacuum module.
As with the case of simple Lie algebras \cite[Ch.~15]{fb:va},
the vertex algebra $C^k(\agot)=V^k(\agot)\ot\Fc$ acquires a structure of a BRST complex
of the quantum Drinfeld--Sokolov reduction. We show that its cohomology
$H^k(\agot)^i$
is zero for all degrees $i\ne 0$ and define the
$\Wc$-algebra by setting $\Wc^{\tss k}(\agot)=H^k(\agot)^0$.

Furthermore, we
give an explicit construction of free generators of the $\Wc$-algebra $\Wc^{\tss k}(\agot)$.
In the particular case $e=0$ they coincide with those previously found in \cite{am:eg}.
A version of the quantum Miura map yields an embedding
of $\Wc^{\tss k}(\agot)$ into the vertex algebra $V^{k+N}(\h)$ associated with the diagonal
subalgebra $\h$ of $\agot$.
In the case $e=0$ we recover the corresponding realization \cite{fl:mt} of $\Wc^{\tss k}(\gl_N)$
as in \cite{am:eg}; see also \cite{a:iw}.

As with this particular case, the vertex algebra $\Wc^{\tss -N}(\agot)$
at the critical level $k=-N$ turns out to be commutative. It is isomorphic
to the classical $\Wc$-algebra $\Wc(\agot)$
introduced in \cite{mr:cw} and to the center
of the vertex algebra $V^{-N}(\agot)$,
as described in \cite{ap:qm} and \cite{m:cc}; cf. \cite{ff:ak}.

Note that if all Jordan blocks of the nilpotent element $e$
are of the same size,
the Lie algebra $\agot$ is isomorphic to a truncated polynomial current algebra
of the form $\gl_n[v]/(v^p=0)$, which is also known as the {\em Takiff algebra}.
This leads to a natural generalization of our definition of the $\Wc$-algebras
to the class of Takiff algebras $\g[v]/(v^p=0)$ associated with
an arbitrary simple Lie algebra $\g$.

\section{BRST cohomology for centralizers}
\label{sec:bc}

Here we adapt the well-known BRST construction of vertex algebras
to the case of centralizers in type $A$. We generally follow
\cite[Sec.~4]{a:rtw} and \cite[Ch.~15]{fb:va} with some straightforward
modifications.

Let $e\in\gl_N$ be a nilpotent matrix and let $\agot$ be
the centralizer of $e$ in $\gl_N$. Suppose that the Jordan canonical form of $e$
has Jordan blocks of sizes
$\la_1,\dots,\la_n$, where
$\la_1\leqslant\dots\leqslant \la_n$ and $\la_1+\dots+\la_n=N$.
The corresponding {\em pyramid} is
a left-justified array of rows of unit boxes such that the top row
contains $\la_1$ boxes, the next row contains $\la_2$ boxes, etc.
Denote by $q^{}_1\geqslant\dots\geqslant q^{}_l$ the column lengths
of the pyramid (with $l=\la_n$).
The {\em row-tableau} is obtained by writing the numbers $1,\dots,N$ into the boxes
of the pyramid consecutively by rows from left to right.
For instance, the row-tableau
\ben
\young(12,345,6789)
\een
corresponds to
the pyramid with the rows of lengths $2,3,4$; its column lengths are
$3,3,2,1$.
We let
$\row(a)$ and $\col(a)$ denote the row and column number of the box containing
the entry $a$.

Denote by
$e_{a\tss b}$ the standard basis elements of the Lie algebra $\gl_N$.
For any $1\leqslant i,j\leqslant n$ and any integral values of $r$ with
$\la_j-\min(\la_i,\la_j)\leqslant r<\la_j$
set
\beql{eijr}
E_{ij}^{(r)}=\sum_{\underset{\scriptstyle\col(b)-\col(a)=r}{\row(a)=i,\ \row(b)=j}} e_{a\tss b},
\eeq
summed over $a,b\in\{1,\dots,N\}$. It is well-known that
the elements $E_{ij}^{(r)}$ form a basis of the Lie algebra $\agot$; see e.g.
\cite{bb:ei} and \cite{ppy:si}.
The commutation relations are given by
\ben
\big[E^{(r)}_{ij},E^{(s)}_{hl}\big]=\de_{hj}\ts E^{(r+s)}_{i\tss l}-\de_{i\tss l}\ts E^{(r+s)}_{hj},
\een
assuming that $E^{(r)}_{ij}=0$ for $r\geqslant \la_j$.

\subsection{Affine vertex algebra}
\label{subsec:vm}

The Lie algebra $\g=\gl_N$ acquires a $\ZZ$-gradation $\g=\bigoplus_{r\in\ZZ} \g_r$
determined by $e$ such that
the degree of the basis element $e_{ab}$ equals $\col(b)-\col(a)$.
We thus get an induced $\ZZ$-gradation
$\agot=\bigoplus_{r\in\ZZ} \agot_r$
on the Lie algebra $\agot$, where $\agot_r=\agot\cap \g_r$. Note that
the element \eqref{eijr} is homogeneous of degree $r$.
The subalgebra $\g_0$ is isomorphic to the direct sum
\beql{godec}
\g_0\cong \gl_{q^{}_1}\oplus\dots\oplus\gl_{q^{}_l}.
\eeq
Equip this subalgebra with the normalized Killing form
\beql{killing}
\langle X,Y\rangle=\frac{1}{2N}\ts\tr\tss(\ad\tss X\ts \ad\tss Y),\qquad X,Y\in \g_0.
\eeq
Now define an invariant symmetric bilinear form on $\agot$
following \cite{ap:qm}.
The value $\langle X,Y\rangle$
for homogeneous elements $X,Y\in \agot$ is found by \eqref{killing} for $X,Y\in\agot_0$, and
is zero otherwise.
Writing $X=X_1+\dots+X_l$ and $Y=Y_1+\dots+Y_l$ in accordance with
the decomposition \eqref{godec}, we get
\ben
\langle X,Y\rangle=\frac{1}{N}\ts\sum_{i=1}^l\big(q_i\ts\tr\tss X_iY_i-\tr\tss X_i\ts\tr\tss Y_i\big).
\een
Therefore, if $\la_i=\la_j$ for some $i\ne j$ then
\ben
\big\langle E_{ij}^{(0)},E_{ji}^{(0)}\big\rangle=
\frac{1}{N}\tss(q_1+\dots+q_{\la_i})=
\frac{1}{N}\ts
\big(\la_1+\dots+\la_{i-1}+(n-i+1)\la_i\big),
\een
and for all $i$ and $j$ we have
\ben
\big\langle E_{ii}^{(0)},E_{jj}^{(0)}\big\rangle=\frac{1}{N}\ts
\big(\de_{ij}(\la_1+\dots+\la_{i-1}+(n-i+1)\la_i)-\min(\la_i,\la_j)\big),
\een
whereas all remaining values of the form on the basis vectors are zero.

The affine Kac--Moody algebra $\wh\agot$
is the central
extension
$
\wh\agot=\agot\tss[t,t^{-1}]\oplus\CC K,
$
where $\agot[t,t^{-1}]$ is the Lie algebra of Laurent
polynomials in $t$ with coefficients in $\agot$. For any $r\in\ZZ$ and $X\in\g$
we will write $X[m]=X\ts t^m$. The commutation relations of the Lie algebra $\wh\agot$
have the form
\ben
\big[X[m],Y[p]\big]=[X,Y][m+p]+m\ts\de_{m,-p}\langle X,Y\rangle\ts K,
\qquad X, Y\in\agot,
\een
and the element $K$ is central in $\wh\agot$.
The vacuum module at the level $k\in\CC$
over $\wh\agot$
is the quotient
\ben
V^k(\agot)=\U(\wh\agot)/\Ir,
\een
where $\Ir$ is the left ideal of $\U(\wh\agot)$ generated by $\agot[t]$
and the element $K-k$.
This module is equipped with a vertex algebra structure
and is known as the ({\em universal}) {\em affine vertex algebra}
associated with $\agot$ and the form $\langle\ ,\ts\rangle$; see \cite{fb:va}, \cite{k:va}.
The vacuum vector is the image of the element $1$ in the quotient
and we will denote it by $\vacu$. Furthermore, introduce the fields
\ben
E^{(r)}_{ij}(z)=\sum_{m\in\ZZ} E^{(r)}_{ij}[m]\ts z^{-m-1}\in \End V^k(\agot)[[z,z^{-1}]]
\een
so that under the state-field correspondence map we have
\ben
Y:E^{(r)}_{ij}[-1]\vacu\mapsto E^{(r)}_{ij}(z).
\een
The map $Y$
extends to the whole of $V^k(\agot)$
with the use of normal ordering.
The translation operator $T$ on $V^k(\agot)$
is determined by the properties
\beql{tra}
T:\vacu\mapsto 0\Fand
\big[T,X[m]\big]=-m\tss X[m - 1],\quad X\in\agot,\quad m<0,
\eeq
where $X[m]$ is understood as the operator of left multiplication by $X[m]$.

\subsection{Affine Clifford algebra}
\label{subsec:aca}

Consider the following triangular decomposition of the Lie algebra $\agot$,
\beql{triang}
\agot=\n_-\oplus\h\oplus \n_+,
\eeq
where the subalgebras are defined by
\ben
\n_-=\text{span of\ }\{E^{(r)}_{ij}\ts|\ts i>j\},\quad
\n_+=\text{span of\ }\{E^{(r)}_{ij}\ts|\ts i<j\}\fand \h=\text{span of\ }\{E^{(r)}_{ii}\},
\een
with the superscript $r$ ranging over all admissible values. Denote by
$\Cl$ the Clifford algebra associated with $\n_+[t,t^{-1}]$, so it is generated
by odd elements $\psi^{(r)}_{ij}[m]$ and $\psi^{(r)*}_{ij}[m]$
with the parameters satisfying the conditions $1\leqslant i<j\leqslant n$ together with
$\la_j-\la_i\leqslant r\leqslant\la_j-1$ and $m\in \ZZ$.
The defining relations
are given by the anti-commutation relations
\ben
\big[\psi^{(r)}_{ij}[m],\psi^{(r)*}_{ij}[-m]\big]=1,
\een
while all other pairs of generators anti-commute. Let $\Fc$ be the Fock representation
of $\Cl$ generated by a vector $\vac$ such that
\ben
\psi^{(r)}_{ij}[m]\tss\vac=0\quad\text{for}\quad m\geqslant 0\Fand
\psi^{(r)*}_{ij}[m]\tss\vac=0\quad\text{for}\quad m> 0.
\een
The space $\Fc$ is a vertex algebra with the vacuum vector $\vac$, and
the translation operator $T$
is determined by the properties $T:\vac\mapsto 0$ and
\ben
\big[T,\psi^{(r)}_{ij}[m]\big]=-m\tss \psi^{(r)}_{ij}[m-1],\qquad
\big[T,\psi^{(r)*}_{ij}[m]\big]=-(m-1)\tss \psi^{(r)*}_{ij}[m-1].
\een
The fields are defined by
\ben
\psi^{(r)}_{ij}(z)=\sum_{m\in\ZZ} \psi^{(r)}_{ij}[m]\ts z^{-m-1}
\Fand
\psi^{(r)*}_{ij}(z)=\sum_{m\in\ZZ} \psi^{(r)*}_{ij}[m]\ts z^{-m}
\een
so that
\ben
Y:\psi^{(r)}_{ij}[-1]\tss\vac\mapsto \psi^{(r)}_{ij}(z)\Fand
Y:\psi^{(r)*}_{ij}[0]\tss\vac\mapsto \psi^{(r)*}_{ij}(z).
\een

The vertex algebra $\Fc$ has a $\ZZ$-gradation $\Fc=\bigoplus_{i\in\ZZ} \Fc^i$,
defined by
\ben
\deg \vac=0,\qquad \deg \psi^{(r)}_{ij}[m]=-1\Fand \deg \psi^{(r)*}_{ij}[m]=1.
\een

\subsection{BRST complex}
\label{subsec:bc}

Introduce the vertex algebra $C^k(\agot)$ as the tensor product
\ben
C^k(\agot)=V^k(\agot)\ot\Fc.
\een
We will use notation $\vacu$ for its vacuum vector $\vacu\ot\vac$.
The vertex algebra $C^k(\agot)$ is $\ZZ$-graded, its $i$-th component has the form
\ben
C^k(\agot)^i=V^k(\agot)\ot\Fc^i.
\een

Consider the fields
$Q(z)$ and $\chi(z)$ defined by
\beql{qz}
Q(z)=\sum_{i<j}E^{(a)}_{ij}(z)\tss \psi^{(a)*}_{ij}(z)
-\sum_{i<j<h} \psi^{(a)*}_{ij}(z) \psi^{(b)*}_{jh}(z) \psi^{(a+b)}_{ih}(z),
\eeq
and
\beql{chiz}
\chi(z)=\sum_{i=1}^{n-1} \psi^{(\la_{i+1}-1)*}_{i\ts i+1}(z).
\eeq
To simplify the formulas,
here and throughout the paper we use the convention that
summation over all admissible values of repeated superscripts of the form $a,b,c$ is assumed.
For instance, summation over $a$ running over
the values $\la_j-\la_i,\dots,\la_j-1$ is assumed within the first sum
in \eqref{qz}. Define the odd endomorphisms $d_{\st}$ and $\chi$ of $C^k(\agot)$
as the residues (coefficients of $z^{-1}$) of the fields \eqref{qz} and \eqref{chiz},
\ben
d_{\st}={\rm res}\ts Q(z)\Fand \chi=\sum_{i=1}^{n-1} \psi^{(\la_{i+1}-1)*}_{i\ts i+1}[1].
\een

\ble\label{lem:diff}
We have the relations
\ben
d_{\st}^2=\chi^2=[d_{\st},\chi]=0.
\een
\ele

\bpf
The relations are verified by the standard OPE calculus with the use
of the Taylor formula and Wick theorem \cite{k:va}.
Using the basic OPEs
\beql{opee}
E^{(r)}_{ij}(z) E^{(s)}_{hl}(w)\sim \frac{1}{z-w}
\Big(\de_{hj}\ts E^{(r+s)}_{i\tss l}(w)-\de_{i\tss l}\ts E^{(r+s)}_{hj}(w)\Big)
+\frac{k\tss\langle E^{(r)}_{ij},E^{(s)}_{hl}\rangle}{(z-w)^2},
\eeq
and
\beql{opepsi}
\psi^{(r)}_{ij}(z)\tss \psi^{(r)*}_{ij}(w)\sim \frac{1}{z-w},\qquad
\psi^{(r)*}_{ij}(z)\tss \psi^{(r)}_{ij}(w)\sim \frac{1}{z-w},
\eeq
we find that the OPE $Q(z)Q(w)$
is regular, thus implying that $d_{\st}^2=0$. The remaining relations are straightforward to
verify.
\epf

By Lemma~\ref{lem:diff}, the odd endomorphism $d=d_{\st}+\chi$
of $C^k(\agot)$
has the properties $d^2=0$ and $d: C^k(\agot)^i\to C^k(\agot)^{i+1}$.
We thus get an analogue $(C^k(\agot)^{\bullet},d)$ of the BRST complex
of the quantum Drinfeld--Sokolov reduction, associated with the Lie algebra $\agot$;
cf. \cite[Ch.~15]{fb:va}.
Since $d$ is the residue of a vertex operator,
the cohomology $H^k(\agot)^{\bullet}$ of the complex
 is a vertex algebra which we will use
to define and describe the $\Wc$-algebras $\Wc^{\tss k}(\agot)$.

\section{$\Wc$-algebras $\Wc^{\tss k}(\agot)$}
\label{sec:dwa}

Introduce another $\ZZ$-gradation on $C^k(\agot)^{\bullet}$ by defining the
(conformal) degrees by
\ben
\deg' E_{ij}^{(r)}[m]=\deg' \psi_{ij}^{(r)}[m]=-m+i-j\Fand
\deg' \psi_{ij}^{(r)*}[m]=-m+j-i.
\een
Observe that the differential $d$ has degree $0$ and so it preserves
this gradation thus defining a $\ZZ$-gradation on the cohomology $H^k(\agot)^{\bullet}$.

\bde\label{def:walg}
The $\ZZ$-graded vertex algebra $H^k(\agot)^0$ is called
the {\em $\Wc$-algebra associated with the centralizer $\agot$
at the level $k$}
and denoted by $\Wc^{\tss k}(\agot)$.
\qed
\ede
Note that the definition also depends on the chosen bilinear form $\langle \ ,\ts \rangle$ on $\agot$.
Our next goal is to prove the following analogue of \cite[Thm~15.1.9]{fb:va}
which describes the structure of principal $\Wc$-algebras associated
with simple Lie algebras.

\bth\label{thm:homol}
The $\Wc$-algebra $\Wc^{\tss k}(\agot)$ is strongly and freely generated by elements
$w_1,\dots,w_N$ of the respective degrees
\ben
\underbrace{1,\dots,1}_{\la_n},\underbrace{2,\dots,2}_{\la_{n-1}},\dots,\underbrace{n,\dots,n}_{\la_1}.
\een
Hence, the Hilbert--Poincar\'e series
of the algebra $\Wc^{\tss k}(\agot)$ is given by
\ben
\prod_{s=0}^{\infty}\ts\prod_{l=1}^n\big(1-q^{\tss l+s}\big)^{-\la_{n-l+1}}.
\een
Moreover, $H^k(\agot)^i=0$ for all $i\ne 0$.
\eth

The proof relies on essentially the same arguments as in
\cite[Ch.~15]{fb:va} (see also \cite[Sec.~4]{a:rtw})
which we will outline in the rest of this section.
A family of generators $w_1,\dots,w_N$ will be produced in
Sec.~\ref{sec:gw}.

For all $1\leqslant i<j\leqslant n$ and $r=\la_j-\la_i,\dots,\la_j-1$
introduce the fields
\beql{eijrz}
e^{(r)}_{ij}(z)=E^{(r)}_{ij}(z)+\sum_{h>j} \psi^{(a)}_{ih}(z) \psi^{(a-r)*}_{jh}(z)
-\sum_{h<i} \psi^{(a)}_{hj}(z) \psi^{(a-r)*}_{hi}(z),
\eeq
where we keep using the convention on the summation over $a$ as in \eqref{qz}.
Similarly, for $i\geqslant j$ and $r=0,1,\dots,\la_j-1$ set
\beql{eijrge}
e^{(r)}_{ij}(z)=E^{(r)}_{ij}(z)+\sum_{h>i} :\psi^{(a)}_{ih}(z) \psi^{(a-r)*}_{jh}(z):
-\sum_{h<j} :\psi^{(a)}_{hj}(z) \psi^{(a-r)*}_{hi}(z):.
\eeq
Note that by the defining relations in the Clifford algebra $\Cl$,
the normal ordering is necessary only for the case where $i=j$ and $r=0$.
Introduce Fourier coefficients $e^{(r)}_{ij}[m]$ of the fields \eqref{eijrz} and \eqref{eijrge}
by setting
\ben
e^{(r)}_{ij}(z)=\sum_{m\in\ZZ} e^{(r)}_{ij}[m]\ts z^{-m-1}.
\een

In the formulas of the next lemmas we assume that the fields with out-of-range
parameters are equal to zero.

\ble\label{lem:commrels}
\begin{enumerate}[label=(\roman*)]
\item
For $i\geqslant j$ and $h<l$ we have
\beql{epsi}
\big[e^{(r)}_{ij}[m],\psi^{(s)*}_{hl}[p]\big]
=\de_{lj}\ts \psi^{(s-r)*}_{h\tss i}[m+p]-\de_{h\tss i}\ts \psi^{(s-r)*}_{jl}[m+p].
\eeq
Moreover, if $i\geqslant j$ and $h\geqslant l$ then
\ben
\big[e^{(r)}_{ij}[m],e^{(s)}_{hl}[p]\big]
=\de_{hj}\ts e^{(r+s)}_{i\tss l}[m+p]-\de_{i\tss l}\ts e^{(r+s)}_{hj}[m+p]
+m\ts\de_{m,-p}\tss(k+N)\tss\langle E^{(r)}_{ij},E^{(s)}_{hl}\rangle.
\een
\item
For $i< j$ and $h<l$ we have
\ben
\big[e^{(r)}_{ij}[m],\psi^{(s)}_{hl}[p]\big]
=\de_{hj}\ts \psi^{(r+s)}_{i\tss l}[m+p]-\de_{i\tss l}\ts \psi^{(r+s)}_{kj}[m+p]
\een
and
\ben
\big[e^{(r)}_{ij}[m],e^{(s)}_{hl}[p]\big]
=\de_{hj}\ts e^{(r+s)}_{i\tss l}[m+p]-\de_{i\tss l}\ts e^{(r+s)}_{hj}[m+p].
\een
\end{enumerate}
\ele

\bpf
All relations are easily verified with the use of the OPEs
\eqref{opee} and \eqref{opepsi}.
\epf

For all $i=1,\dots,n$ set
\beql{alfi}
\al_i=-\la_i+\frac{k+N}{N}\ts\big(\la_1+\dots+\la_{i-1}+(n-i+1)\la_i\big).
\eeq

\ble\label{lem:reldiff}
The following relations hold for all $i\geqslant j$:
\begin{align}
\big[d_{\st},e^{(r)}_{ij}(z)\big]&=\sum_{h=j}^{i-1} :e^{(a+r)}_{hj}(z)\tss \psi^{(a)*}_{hi}(z):
-\sum_{h=j+1}^{i} :\psi^{(a)*}_{jh}(z)\tss e^{(a+r)}_{ih}(z)\ts{:}
+\al_j\tss \de_{r\tss 0}\tss \di_z\tss \psi^{(0)*}_{ji}(z),
\non\\[0.4em]
\big[\chi,e^{(r)}_{ij}(z)\big]&=\psi^{(\la_{i+1}-r-1)*}_{j\ts i+1}(z)
-\psi^{(\la_{j}-r-1)*}_{j-1\ts i}(z).
\label{chieij}
\end{align}
Moreover, for all $i<j$ we have
\begin{alignat}{2}
\big[d_{\st},e^{(r)}_{ij}(z)\big]&=0, &\big[\chi,e^{(r)}_{ij}(z)\big]&=0,
\non\\[0.4em]
\big[d_{\st},\psi^{(r)}_{ij}(z)\big]&=e^{(r)}_{ij}(z),\qquad
&\big[\chi,\psi^{(r)}_{ij}(z)\big]
&=\de_{i\ts j-1}\tss\de_{r\ts\la_{j}-1},
\non
\end{alignat}
and
\ben
\big[d_{\st},\psi^{(r)*}_{ij}(z)\big]=-\sum_{i<h<j}\psi^{(a)*}_{ih}(z)\tss\psi^{(r-a)*}_{hj}(z),
\qquad \big[\chi,\psi^{(r)*}_{ij}(z)\big]=0.
\een
\ele

\bpf
All relations are verified by using the OPEs \eqref{opee} and \eqref{opepsi}.
We give some details for the proof of the first relation. As a first step,
by a direct computation with the use of the Wick theorem we get the OPE
\ben
\bal
Q(z)\tss e^{(r)}_{ij}(w)&\sim \frac{1}{z-w}
\Big(\sum_{h=j}^{i-1} :e^{(a+r)}_{hj}(w)\tss \psi^{(a)*}_{hi}(w):
-\sum_{h=j+1}^{i} :e^{(a+r)}_{ih}(w)\tss\psi^{(a)*}_{jh}(w)\ts{:}\Big)\\
&+\frac{1}{(z-w)^2}\ts\de_{r\tss 0}\ts
\big(k\langle E^{(0)}_{ij},E^{(0)}_{ji}\rangle
+\la_1+\dots+\la_{j-1}+(n-i)\la_j\big)\ts\psi^{(0)*}_{ji}(z),
\eal
\een
where the term $\psi^{(0)*}_{ji}(z)$ is nonzero only if
$j<i$ and $\la_i=\la_j$. Relation \eqref{epsi} of Lemma~\ref{lem:commrels}
implies (assuming summation over $a$) that
\ben
:e^{(a+r)}_{ih}(w)\tss\psi^{(a)*}_{jh}(w)\ts{:}
={:}\ts\psi^{(a)*}_{jh}(w)\tss e^{(a+r)}_{ih}(w)\ts{:}+\de_{r\tss 0}\tss\la_j\tss
\di_w\tss\psi^{(0)*}_{ji}(w).
\een
The required relation now follows by applying
the Taylor formula to $\psi^{(0)*}_{ji}(z)$ to write
\ben
\psi^{(0)*}_{ji}(z)=\psi^{(0)*}_{ji}(w)+(z-w)\tss \di_w\tss\psi^{(0)*}_{ji}(w)+\dots,
\een
and then by taking the residue over $z$ in the resulting expressions.
\epf

Denote by $C^k(\agot)_0$ the subspace of $C^k(\agot)$ spanned by all vectors of the form
\ben
e^{(r_1)}_{i_1j_1}[m_1]\dots e^{(r_q)}_{i_qj_q}[m_q]\tss
\psi^{(s_1)*}_{h_1l_1}[p_1]\dots \psi^{(s_t)*}_{h_tl_t}[p_t]\ts\vacu,\qquad i_a\geqslant j_a,\qquad h_a<l_a,
\een
and by $C^k(\agot)_+$ the subspace of $C^k(\agot)$ spanned by all vectors of the form
\ben
e^{(r_1)}_{i_1j_1}[m_1]\dots e^{(r_q)}_{i_qj_q}[m_q]\tss
\psi^{(s_1)}_{h_1l_1}[p_1]\dots \psi^{(s_t)}_{h_tl_t}[p_t]\ts\vacu,\qquad i_a< j_a,\qquad h_a<l_a.
\een

By Lemma~\ref{lem:commrels}, both $C^k(\agot)_0$ and $C^k(\agot)_+$ are vertex subalgebras
of $C^k(\agot)$. Furthermore, by Lemma~\ref{lem:reldiff} each of the subalgebras
is preserved by the differential $d=d_{\st}+\chi$.
This implies the tensor product decomposition of complexes
\ben
C^k(\agot)^{\bullet}\cong C^k(\agot)^{\bullet}_0\ot C^k(\agot)^{\bullet}_+.
\een
Hence the cohomology of $C^k(\agot)^{\bullet}$ is isomorphic to the tensor product
of the cohomologies of $C^k(\agot)^{\bullet}_0$ and $C^k(\agot)^{\bullet}_+$.

By Lemma~\ref{lem:reldiff}, for $i<j$ we have
\ben
\big[d,e^{(r)}_{ij}[m]\big]=0,\qquad \big[d,\psi^{(r)}_{ij}[m]\big]
=e^{(r)}_{ij}[m]+\de_{i\ts j-1}\tss\de_{r\ts\la_{j}-1}\de_{m,-1}.
\een
Therefore, the complex $C^k(\agot)^{\bullet}_+$ has no higher cohomologies,
while its zeroth cohomology is one-dimensional; see \cite[Sec~15.2.6]{fb:va}.
So the cohomology of $C^k(\agot)^{\bullet}$ is isomorphic to the
cohomology of the complex $C^k(\agot)^{\bullet}_0$. To calculate the latter,
equip this complex with a double gradation by setting
\ben
\bideg e^{(r)}_{ij}[m]=(i-j,j-i),\qquad \bideg \psi^{(r)*}_{ij}[m]=(j-i,i-j+1).
\een
Then $C^k(\agot)^{\bullet}_0$ acquires a structure of bicomplex with
$\bideg \chi=(1,0)$ and $\bideg d_{\st}=(0,1)$. Take $\chi$ as the zeroth
differential of the associated spectral sequence and $d_{\st}$ as the first.
Next we compute the cohomology of $C^k(\agot)^{\bullet}_0$ with respect to $\chi$.

Consider the linear span of all fields $e^{(r)}_{ij}(z)$ with $i\geqslant j$
and $r=0,1,\dots,\la_j-1$.
We will choose a new basis of this vector space which is formed by the fields
\ben
P^{(r)}_l(z)=e^{(r)}_{n\ts n-l+1}(z)+e^{(r+\la_n-\la_2)}_{n-1\ts n-l}(z)
+\dots+e^{(r+\la_n+\dots+\la_{l+1}-\la_{n-l+1}-\dots-\la_2)}_{l\tss 1}(z)
\een
for $l=1,\dots,n$ and $r=0,1,\dots,\la_{n-l+1}-1$ together with
\ben
I^{(r)}_{ij}(z)
=\sum_{h=1}^i e^{(r+\la_{j-1}+\dots+\la_{j-h+1}-\la_i-\dots-\la_{i-h+2})}_{j-h\ts i-h+1}(z)
\een
for $i<j$ and $r=0,1,\dots,\la_i-1$. The following properties of the new basis
vectors are immediate from \eqref{chieij}.

\ble\label{lem:chiact}
We have the relations
\ben
\big[\chi,P^{(r)}_l(z)\big]=0
\Fand
\big[\chi,I^{(r)}_{ij}(z)]=\psi^{(\la_{j}-r-1)*}_{i\tss j}(z).
\vspace{-0.6cm}
\een
\qed
\ele

Lemma~\ref{lem:chiact} allows us to apply the arguments of \cite[Sec.~15.2.9]{fb:va}
to conclude that all higher cohomologies of the complex $C^k(\agot)^{\bullet}_0$
with respect to $\chi$ vanish, while
the zeroth cohomology is the commutative
vertex subalgebra of $C^k(\agot)_0$ spanned by all monomials
\beql{cvs}
P^{(r_1)}_{l_1}[m_1]\dots P^{(r_q)}_{l_q}[m_q]\tss\vacu,
\eeq
where we use the Fourier
coefficients $P^{(r)}_{l}[m]$ defined by
\beql{fc}
P^{(r)}_l(z)=\sum_{m\in\ZZ} P^{(r)}_{l}[m]\tss z^{-m-1}.
\eeq
By a standard procedure outlined in
\cite[Sec.~15.2.11]{fb:va}, each element of this subalgebra gives rise to a unique cocycle
in the complex $C^k(\agot)^{\bullet}_0$ with the differential $d$.
Moreover, the cocycles $W^{(r)}_{l}$ corresponding to the vectors $P^{(r)}_{l}[-1]\tss\vacu$
with $l=1,\dots,n$ and $r=0,1,\dots,\la_{n-l+1}-1$
strongly and freely generate the $\Wc$-algebra $\Wc^{\tss k}(\agot)$.
More precisely, introduce the Fourier
coefficients $W^{(r)}_{l}[m]$ by
\ben
W^{(r)}_l(z)=\sum_{m\in\ZZ} W^{(r)}_{l}[m]\tss z^{-m-1}.
\een
Define a linear ordering on the set of coefficients $W^{(r)}_{l}[m]$ by
positing that the corresponding triples $(l,r,m)$
are ordered lexicographically.
Then the ordered monomials
\beql{basmo}
W^{(r_1)}_{l_1}[m_1]\dots W^{(r_q)}_{l_q}[m_q]\tss\vacu,
\eeq
where all $m_i<0$, form a basis of $\Wc^{\tss k}(\agot)$.
The proof of Theorem~\ref{thm:homol} is completed by the observation that
the conformal degree of the element $W^{(r)}_{l}[m]$ equals $l-m-1$.

\section{Generators of $\Wc^{\tss k}(\agot)$}
\label{sec:gw}

For an $n\times n$ matrix $A=[a_{ij}]$
with entries in a ring we will consider its {\em column-determinant}
defined by
\ben
\cdet\ts A=\sum_{\si\in\Sym_n}\sgn\si\cdot a_{\si(1)\ts 1}\dots a_{\si(n)\ts n}.
\een

We will produce generators of the $\Wc$-algebra $\Wc^{\tss k}(\agot)$ as
elements of the vertex algebra $C^k(\agot)_0$. Combine the Fourier
coefficients $e^{(r)}_{ij}[-1]\in \End C^k(\agot)_0$ into polynomials in
a variable $u$ by setting
\ben
e_{ij}(u)=\sum_{r=0}^{\la_j-1} e^{(r)}_{ij}[-1]\ts u^r,\qquad i\geqslant j.
\een
Let $x$ be another variable and consider the matrix
\ben
\Ec=\begin{bmatrix}x+\al_1\tss T+e_{11}(u)&-u^{\la_2-1}&0&\dots&0\\[0.4em]
                 e_{21}(u)&x+\al_2\tss T+e_{22}(u)&-u^{\la_3-1}&\dots&0\\[0.4em]
                 \dots&\dots&\dots&\dots& \dots \\
                             \dots&\dots&\dots&\dots&-u^{\la_n-1}\\[0.4em]
                             e_{n1}(u)&e_{n2}(u)&\dots&\dots&x+\al_n\tss T+e_{nn}(u)
                \end{bmatrix},
\een
where the constants $\al_i$ are defined in \eqref{alfi}.
Its column-determinant
is a polynomial in $x$ of the form
\beql{cdete}
\cdet\ts\Ec=x^{\tss n}+w_1(u)\tss x^{\tss n-1}+\dots+w_n(u),\qquad w_l(u)=\sum_r w^{(r)}_l\ts u^r,
\eeq
so that the coefficients $w^{(r)}_l$ are endomorphisms of $C^k(\agot)_0$.

The particular case $e=0$ of the following theorem
(that is, with $\la_1=\dots=\la_n=1$) is contained in \cite[Thm~2.1]{am:eg}.

\bth\label{thm:glncent}
All elements $w^{(r)}_l\tss\vacu$ with $l=1,\dots,n$ and
\beql{conda}
\la_{n-l+2}+\dots+\la_n< r+l\leqslant\la_{n-l+1}+\dots+\la_n
\eeq
belong to the $\Wc$-algebra
$\Wc^{\tss k}(\agot)$. Moreover, the $\Wc$-algebra is strongly and freely generated
by these elements.
\eth

\bpf
The first part of the theorem will follow if
we show that the elements $w^{(r)}_l\tss\vacu\in C^k(\agot)_0$
are annihilated by the differential $d$. To verify this property,
it will be convenient to identify $C^k(\agot)_0$ with an isomorphic vertex
algebra $\wt V^k(\agot)$ defined as follows; cf. \cite{am:eg}.
Consider the Lie superalgebra
\beql{affinsu}
\big(\bgot[t,t^{-1}]\oplus \CC K \big)\oplus \m[t,t^{-1}],
\eeq
where the Lie algebra $\bgot$ is spanned by the vectors $e^{(r)}_{ij}$
with $i\geqslant j$ and $r=0,1,\dots,\la_j-1$
understood as basis elements of
the low triangular part $\n_-\oplus\h$
in the decomposition \eqref{triang} via the identification $e^{(r)}_{ij}\rightsquigarrow E^{(r)}_{ij}$,
the even element $K$ is central and $\m$
is the supercommutative Lie superalgebra spanned by (abstract) odd elements $\psi^{(r)*}_{ij}$
with $i<j$ and $r=\la_j-\la_i,\dots,\la_j-1$.
The even component of the Lie superalgebra \eqref{affinsu}
is the Kac--Moody affinization $\bgot[t,t^{-1}]\oplus \CC K$
of $\bgot$ with the commutation relations given by
\ben
\big[e^{(r)}_{ij}[m],e^{(s)}_{hl}[p]\big]
=\de_{hj}\ts e^{(r+s)}_{i\tss l}[m+p]-\de_{i\tss l}\ts e^{(r+s)}_{hj}[m+p]
+m\ts\de_{m,-p}\tss K\tss\langle E^{(r)}_{ij},E^{(s)}_{hl}\rangle,
\een
where the element $e^{(r)}_{ij}[m]$ is now understood as the vector $e^{(r)}_{ij}t^m$.
The remaining commutation relations coincide with those in \eqref{epsi},
where $\psi^{(r)*}_{ij}[m]$ is understood as the vector $\psi^{(r)*}_{ij}t^{m-1}$.
Now define $\wt V^k(\agot)$ as the representation of the Lie superalgebra \eqref{affinsu}
induced from the
one-dimensional representation of
$(\bgot[t]
\oplus \CC K)\oplus \m[t]$
on which
$\bgot[t]$ and $\m[t]$
act trivially and $K$ acts as $k+N$.
Then $\wt V^k(\agot)$ is a vertex algebra isomorphic to $C^k(\agot)_0$
so that the fields with the same names respectively correspond to each other.
Moreover, the cyclic span of the vacuum
vector over the Lie algebra $\bgot[t,t^{-1}]\oplus \CC K$
is a subalgebra of the vertex algebra $\wt V^k(\agot)$ isomorphic to the vacuum module
$V^{k+N}(\bgot)$.

Observe that the coefficients $w^{(r)}_l$ defined in \eqref{cdete}
can now be understood as elements of the universal enveloping algebra
$\U(t^{-1}\bgot[t^{-1}])$. As a vertex algebra, $\wt V^k(\agot)$ is equipped
with the $(-1)$-product, and
each Fourier coefficient $e^{(r)}_{ij}[m]$ with $m<0$ can be regarded as
the operator of left $(-1)$-multiplication by the vector $e^{(r)}_{ij}[m]\tss\vacu$,
and which is the same as the left multiplication by the element $e^{(r)}_{ij}[m]$
in the algebra $\U(t^{-1}\bgot[t^{-1}])$.
Therefore, the monomials in the elements $e^{(r)}_{ij}[m]$ which occur
in the expansion of the column-determinant $\cdet\ts\Ec$ can be regarded as
the corresponding $(-1)$-products calculated consecutively from right to left,
starting from the vacuum vector.

By Lemma~\ref{lem:reldiff}, for $i\geqslant j$ we have the relations
\begin{multline}
\big[d,e^{(r)}_{ij}[-1]\big]=\sum_{h=j}^{i-1} e^{(a+r)}_{hj}[-1]
\psi^{(a)*}_{hi}[0]
-\sum_{h=j+1}^{i} \psi^{(a)*}_{jh}[0] e^{(a+r)}_{ih}[-1]\\[0.3em]
{}+\psi^{(\la_{i+1}-r-1)*}_{j\ts i+1}[0]
-\psi^{(\la_{j}-r-1)*}_{j-1\ts i}[0]+\al_j\tss \de_{r\tss 0}\tss \psi^{(0)*}_{ji}[-1].
\non
\end{multline}
Introducing the Laurent polynomials
\ben
\phi_{ij}=\sum_{r=\la_i-\la_j}^{\la_i-1} \psi^{(r)*}_{ji}[0]\tss u^{-r},\qquad i>j,
\een
we can write the relations in the form
\ben
\big[d,e_{ij}(u)\big]=\Big\{\sum_{h=j}^{i-1} e_{hj}(u)\tss
\phi_{ih}
-\sum_{h=j+1}^{i} \phi_{hj}\tss e_{ih}(u)
+\phi_{i+1\ts j}\tss u^{\la_{i+1}-1}
-\phi_{i\ts j-1}\tss u^{\la_{j}-1}+\al_j\tss T\tss \phi_{ij}\Big\}_+,
\een
where the symbol $\{\dots\}_+$ indicates the component of a Laurent polynomial
containing only nonnegative powers of $u$,
\ben
\big\{\sum_i c_i\tss u^i\big\}_+=\sum_{i\geqslant 0} c_i\tss u^i.
\een

Let $\Ec_{ij}$ denote the $(i,j)$ entry of the matrix $\Ec$. Since $d$ commutes with
the translation operator $T$, we come to the commutation relations
\beql{comdeij}
\big[d,\Ec_{ij}\big]=\Big\{\sum_{h=j}^{i-1} \Ec_{hj}\ts
\phi_{ih}
-\sum_{h=j+1}^{i} \phi_{hj}\ts \Ec_{ih}
+\phi_{i+1\ts j}\tss u^{\la_{i+1}-1}-
\phi_{i\ts j-1}\tss u^{\la_{j}-1}\Big\}_+,
\eeq
which hold for $i\geqslant j$.
The column-determinant of $\Ec$ can be written explicitly in the
form\footnote{This also shows that it coincides with the row-determinant of $\Ec$
defined in a similar way.}
\ben
\cdet\ts\Ec=\sum_{p=0}^{n-1}\ \ \sum_{0=i_0< i_1<\dots<i_p<i_{p+1}=n}
\Ec_{i_1 i_0+1}\ts \Ec_{i_2\ts i_1+1}\dots
\Ec_{i_{p+1}\ts i_p+1}\ts u^{\la_{j_1}-1+\dots+\la_{j_q}-1},
\een
where $\{j_1,\dots,j_q\}$ is the complement to the subset $\{i_0+1,\dots,i_p+1\}$
in the set $\{1,\dots,n\}$.
Since $d$ is the residue of a vertex operator,
$d$ is a derivation of the $(-1)$-product on $\wt V^k(\agot)$.
Hence, using \eqref{comdeij}, we get
\ben
\bal
\big[d,\cdet\ts\Ec\big]&=
\sum_{p=0}^{n-1}\ \ \sum_{0=i_0< i_1<\dots<i_p<i_{p+1}=n}\
\sum_{s=0}^{p}\ts
\Ec_{i_1 i_0+1}\ts \dots \Ec_{i_{s}\ts i_{s-1}+1}\\[0.5em]
{}&\times\Big\{\sum_{i^{}_s<i'_{s+1}<i^{}_{s+1}} \Ec_{i'_{s+1}\ts i^{}_s+1}\ts
\phi^{}_{i^{}_{s+1}\ts i'_{s+1}}
-\sum_{i^{}_s<i'_{s}<i^{}_{s+1}} \phi^{}_{i'_s+1\ts i^{}_s+1}\ts \Ec_{i_{s+1}\ts i'_s+1}\\
&\qquad\qquad\qquad\qquad\qquad\qquad\qquad{}
+\phi^{}_{i_{s+1}+1\ts i_s+1}\tss u^{\la_{i_{s+1}+1}-1}-
\phi^{}_{i_{s+1}\ts i_s}\tss u^{\la_{i_s+1}-1} \Big\}_+\\[0.8em]
&\qquad\qquad\qquad\qquad\qquad\qquad\qquad\quad{}
{}\times\Ec_{i_{s+2}\ts i_{s+1}+1}\dots \Ec_{i_{p+1}\ts i_p+1}\ts u^{\la_{j_1}-1+\dots+\la_{j_q}-1}.
\eal
\een
Now apply
the quasi-associativity property of the $(-1)$-product~\cite[Ch.~4]{k:va},
\ben
(a_{(-1)}b)_{(-1)}c=a_{(-1)}(b_{(-1)}c)
+\sum_{j\geqslant 0}a_{(-j-2)}(b_{(j)}c)+\sum_{j\geqslant 0}b_{(-j-2)}(a_{(j)}c),
\een
to bring the expression to the right-normalized form, where
the consecutive $(-1)$-products are calculated from right to left.
Note that by Lemma~\ref{lem:commrels}\tss(i)
the additional terms coming from
the sums over $j\geqslant 0$ annihilate the vacuum vector because all arising
commutators involve elements with distinct subscripts.

Regarding the above expansion of $[d,\cdet\ts\Ec]$ as written
in the right-normalized form, observe that
if we ignore all symbols $\{\dots\}_+$, then it would turn into a telescoping sum and so
would be identically zero.

As a next step, for a fixed value $l\in\{1,\dots,n\}$ consider the terms
in the expansion of $[d,\cdet\ts\Ec]$ containing the variable $x$ with
the powers at least $n-l$. Such terms can occur only in those summands
where the cardinality of the subset $\{i_0+1,\dots,i_p+1\}$
is at least $n-l+1$. Therefore, the maximum value of the powers
$\la_{j_1}-1+\dots+\la_{j_q}-1$
of the variable $u$ which occur in these terms in the expansion,
equals $\la_{n-l+2}+\dots+\la_n-l+1$. This means that
the coefficients of the powers of $u$ exceeding $\la_{n-l+2}+\dots+\la_n-l$
can be calculated from the expansion $[d,\cdet\ts\Ec]$ with all symbols
$\{\dots\}_+$ omitted. However, as we observed above, this expansion is identically zero.
It is clear from \eqref{cdete} that the degree of the polynomial
$w_l(u)$ equals $\la_{n-l+1}+\dots+\la_n-l$ so that
the relations $d\tss w^{(r)}_l\tss\vacu=0$ hold for the parameters $r$ and $l$
satisfying the conditions of the theorem.

To show that the vectors $w^{(r)}_l\tss\vacu$ are strong and free
generators of $\Wc^{\tss k}(\agot)$,
consider the gradation on
$\U(t^{-1}\bgot[t^{-1}])$ defined by setting the degree of
$e^{(r)}_{ij}[m]$ equal to $j-i$. It is clear from the formulas for
the column-determinant $\cdet\ts\Ec$ that the lowest degree component
of the vector $w^{(r)}_l\tss\vacu$ with $r=r'+\la_{n-l+2}+\dots+\la_n-l+1$
coincides with $P^{(r')}_{l}[-1]\tss\vacu$ for all $r'=0,1,\dots,\la_{n-l+1}-1$, as defined in
\eqref{fc}. Therefore, by the argument completing the proof of Theorem~\ref{thm:homol}
at the end of Sec.~\ref{sec:dwa}, the vector $w^{(r)}_l\tss\vacu$ coincides
with the respective cocycle $W^{(r')}_{l}$.
\epf

\section{Miura map and Fateev--Lukyanov realization}
\label{sec:mm}

Consider the affine Kac--Moody algebra $\wh\h=\h[t,t^{-1}]\oplus\CC\tss K$
associated with $\h$ and the bilinear form defined in Sec.~\ref{subsec:vm}.
Denote its generators by $e^{(r)}_{ii}[m]$ with $i=1,\dots,n$, where
$r$ runs over the set
$0,1,\dots,\la_i-1$ and $m$ runs over $\ZZ$. The element $K$ is central
and the commutation relations
are given by the OPEs
\ben
e^{(r)}_{ii}(z) e^{(s)}_{jj}(w)\sim \frac{K\tss\langle E^{(r)}_{ii},E^{(s)}_{jj}\rangle}{(z-w)^2},
\een
where we set
\ben
e^{(r)}_{ii}(z)=\sum_{m\in\ZZ} e^{(r)}_{ii}[m]\ts z^{-m-1}.
\een
Define the vacuum module $V^{k+N}(\h)$
over the Lie algebra $\wh\h$ as the representation
induced from the
one-dimensional representation of
$\h[t]\oplus \CC K$
on which
$\h[t]$
acts trivially and $K$ acts as $k+N$.
Then $V^{k+N}(\h)$ is a vertex algebra with the vacuum
vector $\vacu$ and translation operator $T$ defined as in \eqref{tra}
for $X\in\h$. Recalling the constants $\al_i$ introduced in \eqref{alfi},
expand the product
\beql{expami}
\big(x+\al_1\tss T+e_{11}(u)\big)\dots
\big(x+\al_n\tss T+e_{nn}(u)\big)=
x^{\tss n}+v_1(u)\tss x^{\tss n-1}+\dots+v_n(u)
\eeq
and define the coefficients $v^{(r)}_l$ by
writing
$
v_l(u)=\sum_r v^{(r)}_l\ts u^r.
$

The particular case $e=0$ of the following proposition
is the realization of the $\Wc$-algebra $\Wc^{\tss k}(\gl_n)$
given by Fateev and Lukyanov~\cite{fl:mt}; see also \cite{am:eg}.

\bpr\label{prop:ff}
The elements $v^{(r)}_l\tss\vacu$ with $l=1,\dots,n$ and
$r$ satisfying \eqref{conda}
generate a subalgebra of the vertex algebra $V^{k+N}(\h)$,
isomorphic to the $\Wc$-algebra
$\Wc^{\tss k}(\agot)$.
\epr

\bpf
The Lie algebra projection
$\bgot \to\h$ with the kernel $\n_-$
induces the vertex algebra
homomorphism
$V^{k+N}(\bgot)\to V^{k+N}(\h)$. As we have seen in the previous section,
the $\Wc$-algebra $\Wc^{\tss k}(\agot)$ can be regarded as a subalgebra
of the vertex algebra $V^{k+N}(\bgot)$.
Hence, we get a vertex algebra
homomorphism
\ben
\Upsilon:\Wc^{\tss k}(\agot)\to   V^{k+N}(\h),
\een
obtained by restriction,
which we can call the {\em Miura map}; cf. \cite[Sec.~5.9]{a:iw}, \cite[Sec.~15.4]{fb:va}.
For the image of the column-determinant \eqref{cdete} we have
\ben
\Upsilon:\cdet\ts\Ec\mapsto \big(x+\al_1\tss T+e_{11}(u)\big)\dots
\big(x+\al_n\tss T+e_{nn}(u)\big).
\een
Therefore, the Miura images of the generators of $\Wc^{\tss k}(\agot)$
provided by Theorem~\ref{thm:glncent}
are found by
\ben
\Upsilon: w^{(r)}_l\tss\vacu\mapsto v^{(r)}_l\tss\vacu.
\een

It remains to verify that the Miura map is injective.
For $l=1,\dots,n$ and
$r$ satisfying \eqref{conda}
introduce the Fourier
coefficients $w^{(r)}_{l}[m]$ by
\ben
w^{(r)}_l(z)=\sum_{m\in\ZZ} w^{(r)}_{l}[m]\tss z^{-m-1}.
\een
By Theorem~\ref{thm:glncent},
the monomials
\beql{basmogen}
w^{(r_1)}_{l_1}[m_1]\dots w^{(r_q)}_{l_q}[m_q]\tss\vacu,\qquad m_i<0,
\eeq
with the factors ordered as in \eqref{basmo},
form a basis of $\Wc^{\tss k}(\agot)$.
Suppose that a certain linear combination $w$ of these monomials
belongs to the kernel of $\Upsilon$. That is,
the corresponding linear combination
of their images is zero in the vertex algebra $V^{k+N}(\h)$.
Then so is the linear combination of the components
of the images of the top degree, regarded as elements of the symmetric algebra
$\Sr(t^{-1}\h[t^{-1}])\cong \gr V^{k+N}(\h)$, where we equip
the vacuum module $V^{k+N}(\h)$ with
the canonical filtration of the universal enveloping algebra.
However,
it follows from the proof of \cite[Prop.~4.3]{mr:cw} that
the top degree components of the elements
$T^s\tss v^{(r)}_l$, where
$s\geqslant 0$ and $l=1,\dots,n$ with $r$ satisfying conditions \eqref{conda},
are algebraically independent.
This implies that the linear combination $w$ is zero.
\epf

At the critical level $k=-N$ we have $\al_i=-\la_i$ for all $i=1,\dots,n$ so that
the change of signs $e_{ii}(u)\mapsto -e_{ii}(u)$ allows up to identify
the coefficients in the expansion \eqref{expami} with the generators
of the classical $\Wc$-algebra $\Wc(\agot)$ and with the generators of the
center $\z(\wh\agot)$ of the vertex algebra $V^{-N}(\agot)$ found in
\cite{m:cc} and \cite[Sec.~4]{mr:cw}.

\bco\label{cor:clce}
The $\Wc$-algebra $\Wc^{-N}(\agot)$ is a commutative vertex algebra and
we have isomorphisms
\ben
\Wc^{-N}(\agot)\cong \z(\wh\agot)\cong \Wc(\agot).
\vspace{-0.6cm}
\een
\qed
\eco
In the case $e=0$ we recover the corresponding result of \cite{ff:ak}
in type $A$.

%\newpage
\bigskip\bigskip

\small

\noindent
School of Mathematics and Statistics\newline
University of Sydney,
NSW 2006, Australia\newline
alexander.molev@sydney.edu.au

\end{document}